\documentclass[11pt, a4paper, oneside, leqno]{article}

\usepackage[utf8]{inputenc}

\usepackage{geometry}
\usepackage{fancyhdr}
\usepackage{enumerate}

\usepackage{amsthm}
\usepackage{amsmath}
\usepackage{mathtools}
\usepackage{amssymb}
\usepackage{latexsym}

\theoremstyle{definition}

\theoremstyle{plain}

\theoremstyle{plain}

\usepackage{reversemath}

\usepackage[round,elide]{natbib}
\bibpunct{(}{)}{,}{a}{}{,}

\usepackage[final]{pdfpages}
\usepackage[bookmarks,colorlinks,breaklinks]{hyperref}
\usepackage{doi}

\newcommand{\gcode}[1]{\ulcorner{}#1\urcorner{}}

\hypersetup{
    linkcolor=blue,
    citecolor=blue,
    filecolor=blue,
    urlcolor=blue
}

\title{From foundations to applications: \\ reverse mathematics and philosophy}
\author{%
    Benedict Eastaugh%
    \thanks{%
        University of Warwick,
        \href{mailto:benedict@eastaugh.net}{benedict@eastaugh.net}.
        Thanks to Marianna Antonutti Marfori, Walter Dean, and Alex Paseau for
        their helpful comments.
        This work was supported by the Arts \& Humanities Research Council
        [grant number UKRI3466].
    }
}
\date{July 7, 2026}

\begin{document}


\maketitle

Reverse mathematics is a branch of mathematical logic dedicated to determining
the minimal set existence principles necessary and sufficient to derive ordinary
mathematical theorems about concrete structures like the real line. Since the
mid-1970s, reverse mathematics has developed a systematic classification of the
strength of theorems in areas of mathematics ranging from real and complex
analysis to infinitary combinatorics.
This essay will place reverse mathematics in its historical and philosophical
context, and reveal its relevance to central issues in the philosophy of
mathematics, from the foundational programmes of Hilbert and Brouwer to
contemporary debates about realism, determinacy, and applicability of
mathematics.
In doing so, it will discuss the role of computability theory in measuring the
strength of set existence principles, as well as related questions about
idealisation when these principles are applied in the physical sciences and in
philosophy.


\section{From theorems to axioms}
\label{sec:history}

Establishing equivalences between axioms and theorems is not a new activity. We
can trace it back to the study of Euclid's parallel postulate and the proposals
for substitute principles such as Playfair's axiom, as already noted in late
antiquity by Proclus.%
\footnote{%
\citet{Heath1908,Lewis1920}.
}
The development of mathematical logic in the late 19th and early 20th centuries
made it possible to understand such equivalences as provable biconditionals of
the form
\begin{equation*}
    B \vdash \varphi \leftrightarrow \psi
\end{equation*}
where $\varphi$ and $\psi$ are sentences in the language of the base theory $B$
and $\vdash$ denotes the provability relation of a logic such as classical
first-order logic. The non-trivial cases of such equivalences are those in which
$\varphi$ (and hence also $\psi$) are independent of $B$,
i.e.\ $B \not\vdash \varphi$ and $B \not\vdash \neg\varphi$.
Many readers will be familiar with one such equivalence, between the axiom of
choice and the well-ordering principle.%
\footnote{%
\citet{Zermelo1904,Zermelo1908}.
}

Reverse mathematics is a subfield of mathematical logic which investigates these
equivalences within the formal setting of \emph{second-order arithmetic}.
Second-order arithmetic is a quantified two-sorted language with the
non-logical symbols of first-order arithmetic plus the membership symbol $\in$.
This links the two sorts by allowing us to express that a certain number (given
by a term of the first sort) is a member of a certain set (given by a term of
the second sort). The proof theory and semantics of second-order arithmetic are
those of two-sorted first-order logic, i.e.\ Henkin semantics.
The main objects of study in reverse mathematics are \emph{subsystems of
second-order arithmetic}. These are axiom systems formulated in $\lang_2$ all of
whose axioms are theorems of the theory $\mathsf{Z}_2$ of second-order
arithmetic, a second-order formalisation of the Dedekind--Peano axioms including
the comprehension scheme for all formulas in the second-order language.

Much research in reverse mathematics has been dominated by the `Big Five'
subsystems of second-order arithmetic: $\RCA_0$, $\WKL_0$, $\ACA_0$, $\ATR_0$,
and $\PiCA_0$. These systems were introduced by \citet{Friedman1975,
Friedman1976}, who also proved a number of equivalences between theorems of
analysis and combinatorics, and the characteristic axioms of these subsystems.%
\footnote{%
The historical development of the subsystems is related by
\citet{DeanWalsh2017}.
}
Although its roots lie deeper, in the arithmetisation of analysis by Dedekind
and Cantor, second-order arithmetic as a formal system and its use for
formalising analysis first emerge in the second volume of Hilbert and Bernays's
\emph{Grundlagen der Mathematik} \citeyearpar{HilbertBernays1939}.%
\footnote{%
Supplement~IV, `Formalismen zur deduktiven Entwicklung der Analysis',
pp.~467--512 of \citep{HilbertBernays1970}.
}

Second-order arithmetic is thus, even in its historical inception, closely
connected not only with the broader issue of foundations of mathematics, but
with the Hilbert programme and finitism in particular.
It is a highly expressive language and can formalise much of countable and
countably representable mathematics, including many classical theorems of real,
complex, and functional analysis; theorems from algebra such as the theories of
countable commutative rings, abelian groups, and ordered fields; and countably
infinite combinatorics such as Ramsey's theorem and König's infinity lemma.
Reverse mathematics studies the equivalences between these theorems and axioms
asserting the existence of different classes of infinite, non-computable sets,
such as the axiom schemes of arithmetical comprehension or arithmetical
transfinite recursion.

Results in reverse mathematics are thus closely connected to computability
theory and the hierarchies of sets it studies, such as the arithmetical and
hyperarithmetical hierarchies. The base theory $\RCA_0$, in which the
equivalences of reverse mathematics are typically proved, also has a
relationship to computability, as its characteristic axiom asserts the existence
of all computable sets.
Computability provides a connection to constructive mathematics, as the
mathematical theory of computation has been used to make precise the notion of
construction. This idea also provides a way of demonstrating that a classical
mathematical theorem is not constructively provable, by exhibiting what is known
as a \emph{recursive counterexample}. This is a standard technique in computable
and constructive analysis, and because recursive counterexamples are computable
objects, their existence can be proved in the base theory $\RCA_0$, forming the
basis of many reversals from theorems to axioms.

Recursive counterexamples were developed from the 1940s onwards, while the proof
theory of second-order arithmetic and its subsystems developed in parallel. The
two areas were brought together decisively to form the new field of reverse
mathematics in the mid-1970s. Harvey Friedman introduced the idea of reverse
mathematics in a talk at the International Congress of Mathematicians in
1974 \citep{Friedman1975}. It developed into a substantial research area in
mathematical logic in the 1980s and 1990s, particularly through the work of
Stephen Simpson and his students and collaborators. In recent years the field
has become more diverse in its methodology and its outlook, encompassing work on
computability-theoretic reducibility notions and higher-order reverse
mathematics.%
\footnote{%
These developments are surveyed in \citet{DzhafarovMummert2022}
and \citet{Eastaugh2024}.
}
Because it has received the most attention from philosophers, we will
concentrate on `classical' reverse mathematics, meaning provable equivalences
between axioms and theorems over the base theory $\RCA_0$.


\section{Subsystems of second-order arithmetic}
\label{sec:subsystems}

Second-order arithmetic $\lang_2$ is a two-sorted quantificational language.
Variables $x_0,x_1,\dots$ of the first sort are called \emph{number variables}
and, in the intended interpretation, range over natural numbers.
Variables $X_0,X_1,\dots$ of the second sort are called \emph{set variables}
since in their intended interpretation they range over sets of natural numbers.
The non-logical vocabulary consists of the symbols familiar from first-order
Peano arithmetic: constants $0$ and $1$, function symbols $+$ and $\times$, and
the less-than relation symbol $<$, plus the membership relation symbol $\in$.
The \emph{numerical terms} consist of the number variables, the constants $0$
and $1$, and any term of the form $t_1 + t_2$ or $t_1 \times t_2$ where $t_1$
and $t_2$ are numerical terms. The \emph{atomic formulas} of $\lang_2$ are all
expressions of the form $t_1 = t_2$, $t_1 < t_2$, and $t_1 \in X$, where
$t_1$ and $t_2$ are numerical terms and $X$ is a set variable. Notice that there
are no atomic formulas corresponding to identity for set variables. This is
instead defined in terms of co-extensionality,
\begin{equation*}
    \label{eq:set_var_extensionality}
    \tag{$=_1$}
    X = Y \Leftrightarrow \forall{x}(x \in X \leftrightarrow x \in Y).
\end{equation*}

The \emph{formulas} of $\lang_2$ are obtained by closing the atomic formulas
under propositional connectives and universal and existential quantifiers. When
a quantifier binds a number variable we call it a \emph{number quantifier},
while when it binds a set variable we call it a \emph{set quantifier}. Formulas
of second-order arithmetic in prenex normal form are stratified into a
hierarchy. A formula is $\Sigma^0_0$ if it contains no set quantifiers and all
occurrences of number quantifiers are bounded, i.e.\ have the form
$\exists{x}(x < t \wedge \varphi(x))$ or
$\forall{x}(x < t \rightarrow \varphi(x))$. A prenex formula is $\Sigma^0_{n+1}$
if it has the form $\exists{x}\varphi(x)$ where $\varphi$ is $\Pi^0_n$, while it
is $\Pi^0_{n+1}$ if it has the form $\forall{x}\psi(x)$ where $\psi$ is
$\Sigma^0_n$. $\Sigma^0_n$ and $\Pi^0_n$ formulas are called \emph{arithmetical}
since they do not involve quantification over sets of numbers, although they may
include free set variables.

The semantics of second-order arithmetic are first-order, in the following
sense. An $\lang_2$-structure has the form
\begin{equation*}
    \mathcal{M} = \str{
        M, \mathcal{S}^\mathcal{M},
        0^\mathcal{M}, 1^\mathcal{M},
        +^\mathcal{M}, \times^\mathcal{M},
        <^\mathcal{M}
    }
\end{equation*}
where $M$ is a non-empty set over which the number variables range and
$\mathcal{S}^\mathcal{M} \subseteq \mathcal{P}(M)$ is a non-empty set of subsets
of $M$ over which the set variables range.%
\footnote{%
In other words, general or Henkin semantics rather than the standard semantics
for second-order logic, in which second-order variables always range over the
entire powerset $\mathcal{P}(M)$ of the first-order domain
\citep[\S\S 4.2--4.3]{Shapiro1991}.
}

The \emph{base theory} in which equivalences between (formalisations of)
mathematical theorems and axioms are proved plays a key role in reverse
mathematics. It must be weak enough that the theorems and axioms in question are
not already provable, but strong enough to carry out the proof of the
equivalence. The standard base theory used in reverse mathematical practice is
called $\RCA_0$. Its axioms include, as do all the systems discussed here, the
basic arithmetical axioms of $\PA^-$, i.e.\ those of Peano arithmetic minus the
first-order induction scheme.%
\footnote{%
See \citet[p.~4]{Simpson2009}.
}
$\RCA_0$ also includes the $\Sigma^0_1$ induction scheme, i.e.\ all universal
generalisations of formulas of the form
\begin{equation*}
    \tag{$\Sigma^0_1\text{-}\mathsf{IND}$}
    \varphi(0) \wedge \forall{n}(\varphi(n) \rightarrow \varphi(n+1))
    \rightarrow
    \forall{n}\varphi(n)
\end{equation*}
where $\varphi(n)$ is a $\Sigma^0_1$ formula, possibly with additional free
number and set variables.
Finally, $\RCA_0$ includes the \emph{recursive comprehension axiom} scheme
which gives the system its name, i.e.\ all universal generalisations of formulas
of the form
\begin{equation*}
    \tag{$\Delta^0_1\text{-}\mathsf{CA}$}
    \forall{n}(\varphi(n) \leftrightarrow \psi(n))
    \rightarrow
    \exists{X}\forall{n}(n \in X \leftrightarrow \varphi(n))
\end{equation*}
where $\varphi(n)$ is a $\Sigma^0_1$ formula and $\psi(n)$ is a $\Pi^0_1$
formula, which like the formulas in the induction scheme may have free formula
and set variables, except that the set variable $X$ may not occur freely in
them.

The method of arithmetisation, familiar from its roots in the work of Cantor and
Dedekind, makes clear how to formalise statements about natural, rational, and
real numbers in the language of second-order arithmetic.%
\footnote{%
\citet{Cantor1872, Dedekind1872}. \citet[p.~37 ff.]{Dauben1979} gives a
condensed presentation of Cantor's theory of irrational numbers, while
\citet[p.~29 ff.]{Hallett1984} provides a critical assessment of
arithmetisation. 
For more precise details on the coding operations described here, see
\S\S II.2--II.4 of \citet{Simpson2009}, or \S\S 3.2--3.3 of
\citet{Eastaugh2024}.
}
An arithmetically definable pairing function such as $(m,n) = (m + n)^2 + m$
allows pairs of numbers to be coded by single numbers. This simple operation is
the foundation of all else. Integers can be represented as pairs of natural
numbers, together with operations of integer addition, subtraction, and
multiplication, which treat each pair $(m,n)$ as the sum $m - n$. Rational
numbers are pairs of integers $(n,d)$ with $d \neq 0$, interpreted as the
quotient $\frac{n}{d}$. Real numbers are coded as infinite sequences of
rational numbers obeying a version of the Cauchy convergence criterion, where an
infinite sequence is a (code for) a function $f : \N \to \Q$. Functions are
represented by sets of (codes of) pairs $(m,n)$ such that for every $m$ there
exists exactly one $n$ such that $(m,n) \in f$.

Using these representations, many basic facts about the ring of integers and
the rational and real number fields can be proved, for example that $\Q$ is an
ordered field, or that $\R$ is archimedean. More substantial facts can also be
proved in $\RCA_0$, such as a form of completeness of the real numbers called
nested interval completeness, or a form of the Baire category theorem for
$\R^n$. The uncountability of the real numbers, in the form of the statement
``For every countably infinite sequence of real numbers, there exists a real
number which does not occur in the sequence'', is provable in $\RCA_0$ by
formalising the standard diagonalisation, since the diagonal construction itself
is computable relative to the initial sequence of reals.

$\RCA_0$ is closely linked to computable mathematics. Its characteristic axiom,
the axiom scheme of recursive comprehension, is so called because of Post's
theorem that the sets of natural numbers which are $\Delta^0_1$ definable
(meaning those definable by both a $\Sigma^0_1$ and a $\Pi^0_1$ formula) are
exactly the computable (historically called \emph{recursive}) sets of natural
numbers. The axioms of $\RCA_0$ are \emph{computably true}, meaning that they
are true when the set quantifiers in the axioms are interpreted as ranging over
the computable sets, with the number quantifiers ranging over the standard
natural numbers $\omega = \left\{ 0,1,2,\dots \right\}$ and the symbols $0$,
$1$, $+$, $\times$, and $<$ having their standard meanings. Structures of this
sort, with a standard first-order part, are called \emph{$\omega$-models}.
Since they are distinguished from one another entirely by their second-order
parts, they are usually referred to only in terms of that second-order part:
the $\omega$-model $\REC$ whose second-order part is the class $\REC$ of
computable sets, the $\omega$-model $\ARITH$ whose second-order part is the
class $\ARITH$ of arithmetically definable sets, and so on.

The system $\WKL_0$ is obtained by adding to the axioms of $\RCA_0$ an
additional axiom known as \emph{weak König's lemma} ($\mathrm{WKL}$). This
axiom is a restriction to countably infinite trees formed of sequences of $1$s
and $0$s of König's lemma, the combinatorial principle familiar from graph
theory and set theory which states that every finitely branching infinite tree
has an infinite path through it. To see how to state this in the language of
second-order arithmetic, we start by noting that in $\RCA_0$ one can code finite
sequences of natural numbers by individual numbers.%
\footnote{%
This can be done by, for example, using Gödel's $\beta$-function which may be
familiar from proofs of the incompleteness theorems.
}
$2^{<\N}$ denotes the set of all finite sequences of $0$s and $1$s.
A set $T \subseteq 2^{<\N}$ is a \emph{tree} if for all $t \in T$, if $s$ is an
initial subsequence of $t$, then $s \in T$. A function $f : \N \to 2^{<\N}$ is a \emph{path} through $T$.
Weak König's lemma or $\mathrm{WKL}$ is the $\lang_2$-statement that every
infinite tree $T \subseteq 2^{<\N}$ has a path.

$\WKL_0$ can prove versions of the Heine--Borel theorem for sequential covers,
the Hahn--Banach theorem for separable spaces, as well as key theorems from
mathematical logic such as Gödel's completeness theorem and the compactness
theorem for propositional and first-order logic. All of these theorems are
equivalent to weak König's lemma over $\RCA_0$, showing that they are not
derivable in $\RCA_0$.
This is because there exist models of $\RCA_0$ that are not models of $\WKL_0$,
most notably the $\omega$-model $\REC$. This follows from a construction due to
\citet{Kleene1952a} of what has come to be called the \emph{Kleene tree}, a
computable tree $T \subseteq 2^{<\N}$ with no computable path.%
\footnote{%
For the standard construction of a Kleene tree, see theorem~9.3.2 of
\citet{Soare2016}. It relies on the existence of computably inseparable sets,
as introduced by \citet{Kleene1950}.
}
$\WKL_0$ includes recursive comprehension, so it can prove that Kleene trees
exist, but since the paths through such trees are not computable, weak König's
lemma is false in $\REC$.

The \emph{arithmetical comprehension axiom scheme} consists of the universal
closures of all formulas of the form
\begin{equation*}
    \tag{$\mathrm{ACA}$}
    \exists{X}\forall{n}(n \in X \leftrightarrow \varphi(n))
\end{equation*}
where $\varphi$ is an arithmetical formula and $X$ is not free in $\varphi$.
$\ACA_0$ is the subsystem of second-order arithmetic whose axioms are those of
$\RCA_0$ plus all instances of the arithmetical comprehension axiom scheme.
This system is substantially stronger than both $\RCA_0$ and $\WKL_0$.
It can prove the sequential completeness and compactness of the real numbers:
the Bolzano--Weierstra\ss{} theorem, the monotone convergence theorem, and
generalisations of these theorems to arbitrary complete separable metric spaces.

$\ACA_0$ can also prove König's infinity lemma, the statement that any finitely
branching infinite tree on $\N$ has an infinite path. `Full' König's lemma is
strictly stronger than weak König's lemma. This is because there are
computable, finitely branching infinite trees $T$ on $\N$ such that every path
$p$ through $T$ computes the \emph{halting problem}
\begin{equation*}
    K = \left\{ e : \Phi_e(e)\conv \right\},
\end{equation*}
i.e.\ the set of indexes of Turing machines which halt on every input.%
\footnote{%
$\Phi_e$ denotes the function computed by the Turing machine with index $e$. The
notation $\Phi_e(n)\conv$ means that the Turing machine with index $e$ halts
when given the input $n$.
}
Since there are $\omega$-models of $\WKL_0$ which do not contain $K$, $\WKL_0$
does not prove König's infinity lemma.%
\footnote{%
The existence of such $\omega$-models of $\WKL_0$ follows from the low basis
theorem of \citet{JockuschSoare1972}. For textbook presentations see
\citet[\S 3.7.2 and chapter~9]{Soare2016},
\citet[p.~59]{Hirschfeldt2014},
and \citet[\S VIII.2]{Simpson2009}.
}
On the other hand, arithmetical comprehension is equivalent over $\RCA_0$ to
König's lemma. As a consequence, every $\omega$-model of $\ACA_0$ must contain
$K$ and indeed all of the finite iterations of the \emph{Turing jump} operation
\begin{equation*}
    \jump{X} = \left\{ e : \Phi_e^X(e)\conv \right\},
\end{equation*}
the relativisation of the halting problem to an oracle $X$. The minimal
$\omega$-model of $\ACA_0$ has the class of arithmetically definable sets
$\ARITH$ as its second-order part; by Post's theorem these are exactly the sets
which are computable relative to some finite number of iterations of the Turing
jump to the empty set.

We will not discuss the fourth and fifth systems of the Big Five in great
detail, since many of the philosophical issues connected to reverse mathematics
already appear in relation to the first three systems, but for the sake of
completeness and the few times we will need to discuss them, we briefly
introduce them now.
The fourth system of the Big Five is $\ATR_0$. Its characteristic axiom is the
scheme of \emph{arithmetical transfinite recursion}. This has a slightly
technical definition, but roughly speaking it is the principle that if a set
$X$ codes a well-ordering, then one can iterate any arithmetical operation along
that well-ordering.%
\footnote{%
For a precise definition, see \citet[\S V.2]{Simpson2009}
}
$\ATR_0$ is thus strictly stronger than $\ACA_0$, since not only do its axioms
extend arithmetical definability into the transfinite, it can prove the
existence of the minimal $\omega$-model $\ARITH$ of $\ACA_0$.

The theorems provable in $\ATR_0$ extend beyond real and complex analysis and
into the lower reaches of descriptive set theory. They include Lusin's
separation theorem, that any two disjoint analytic sets can be separated by a
Borel set; the perfect set theorem, that every uncountable closed set has a
perfect subset; and the fact that any two countable well-orderings are
comparable.
$\PiCA_0$ is the final member of the Big Five. Its characteristic axiom is the
$\Pi^1_1$ comprehension axiom. This asserts the existence of sets definable by
$\Pi^1_1$ formulas, those of the form $\forall{X}\varphi(X,n)$ where $X$ is a
set variable and $\varphi$ is an arithmetical formula.
$\PiCA_0$ is strictly stronger than $\ATR_0$. The theorems it can prove include
more results from descriptive set theory, such as the Cantor--Bendixson theorem.


\section{Constructivity and computability}
\label{sec:constructivism}

Much of the early philosophical discussion of reverse mathematics revolved
around its relevance for the foundations of mathematics, and in particular its
connections to the foundational programmes of intuitionism, predicativism, and
finitism developed in the first decades of the 20th century. Despite the
substantial intellectual differences between these programmes, they nevertheless
share one important characteristic, namely an adherence to a spirit of what can
broadly be called constructivism: the view that the only mathematical objects
which exist are those which can be constructed.
There is at least a surface-level similarity with the programme of reverse
mathematics, since the set existence principles that characterise different
subsystems of second-order arithmetic are computability-theoretic or
definability-theoretic in nature. Perhaps unsurprisingly, there are in fact
deeper connections between foundational perspectives and particular systems.
However, these relationships are not straightforward: one cannot simply read
off the mathematical resources of a given foundational standpoint by looking at
the theorems provable in a potentially associated subsystem of second-order
arithmetic.

The first connection is that between computable functions and the various forms
of constructivism, in the more restricted sense particular to Brouwer's
intuitionism, the Russian school of constructive mathematics following Markov,
and Bishop's constructive mathematics. From the 1940s onwards, researchers in
recursive function theory, as computability theory was then called, attempted to
use results concerning computable functions in order to better understand the
limits of constructivism as a foundation.

One difficulty with the identification of constructive with computable is that
constructive mathematicians do not accept the law of the excluded middle, and so
not everything that is computably true (that is, true in the $\omega$-model
$\REC$) is constructively provable. A standard example of this is the
intermediate value theorem, which is true in $\REC$ and provable in $\RCA_0$,%
\footnote{%
\citet[theorem~II.6.6]{Simpson2009}.
}
but which is nonetheless non-constructive, since its proof uses a
non-constructive case distinction in an essential way. Despite the seeming
naturalness of identifying constructivity with computability, examples like this
one show that computable truth is not sufficient for constructive provability.

Depending on the version of constructivism in question, it may not be necessary
either. Brouwer's conception of construction included not just \emph{lawlike}
potentially infinite sequences or functions given by rules, but also
\emph{absolutely free} or \emph{lawless} choice sequences in which the
intuitionistic mathematician can freely choose a value at every step.%
\footnote{%
For an introduction to free choice sequences, see chapter~3 of
\citep{vanAtten2004}.
}
When interpreted in a classical model, such sequences can be non-computable, as
seen in Kleene's \citeyearpar{Kleene1952a} work on Brouwer's fan theorem. The
contrapositive of the fan theorem is weak König's lemma, and thus when the fan
theorem is included in a classical system, it implies the existence of
non-computable sets.
Other forms of constructivism aim to be compatible with \emph{constructive
Church's thesis} (CT), the claim that every function $f : \N^k \to \N$ is
computable. As well as the recursive constructive mathematics of Markov, this
includes Bishop's constructive mathematics, despite Bishop's opposition to the
identification of constructive mathematics with computable mathematics.%
\footnote{%
Representative remarks appear in \citep[pp.~6,~74]{Bishop1967} and
\citep[p.~514]{Bishop1975}.
}

Cardinality considerations show us that most real numbers are non-computable,
since there are uncountably many reals and only countably many computable
reals.%
\footnote{%
A computable real is, roughly speaking, a real number whose value can be
calculated by a computer program to any desired degree of precision. For details
see \citet[p.~14]{Pour-ElRichards1989}.
}
In other words, there are many `gaps' in the computable real numbers $\R^\REC$.
One way in which this manifests itself is in the fact that many theorems of
classical analysis are false when the range of the quantifiers in the statements
of these theorems are restricted to the computable reals. Consider a theorem
$\theta$ of the form
\begin{equation*}
    \forall{X}(\varphi(X) \rightarrow \exists{Y}\psi(X,Y)).
\end{equation*}
A \emph{recursive counterexample} to $\theta$ is a computable set which
satisfies its antecedent but not its consequent. In other words, it is a
computable set $X \subseteq \omega$ such that $\varphi(X)$ holds but there is
no computable $Y \subseteq \omega$ such that $\psi(X,Y)$ holds. The existence
of a recursive counterexample to a theorem shows that it is computably false.

For constructivists who take CT to be consistent, non-computability implies
non-constructivity. To show that a statement is constructively unprovable, it
therefore suffices to show that there exists a recursive counterexample to it.
An early example was Specker's \citeyearpar{Specker1949} construction of a
computable, bounded, monotone sequence of rational numbers whose limit is
non-computable. The existence of \emph{Specker sequences}, as they are now
known, shows that the monotone convergence theorem---together with other
theorems which express the sequential completeness of the real line---is
computably false, and thus constructively unprovable.
In reverse mathematics, recursive counterexamples also play another role:
reversals are often mediated by the construction of a recursive counterexample
in the base theory. For example, one can construct a Specker sequence using
recursive comprehension, apply the monotone convergence theorem to obtain its
limit, which must compute the set $K$ of solutions to the halting problem,
thereby implying $\Sigma^0_1$ comprehension and hence arithmetical
comprehension.%
\footnote{%
For details of these arguments, see lemma~III.1.3 and theorem~III.2.2 in
\citet[pp.~105--107]{Simpson2009}.
}

Another example is a version of the Heine--Borel theorem ($\mathrm{HB}$) stating
that every countable open cover of the closed unit interval has a finite
subcover. A recursive counterexample to $\mathrm{HB}$ is known as a
\emph{singular cover}, and is constructed by forming a computable sequence of
open intervals $C_n = (l_n,r_n)$ which covers every computable real number, but
which has classical measure $<1$. There will then be reals $x \in [0,1]$ not
covered by any $C_n$, and which therefore must be non-computable
\citep{KreiselLacombe1957,LeRouxZiegler2008}. The construction of this recursive
counterexample can be formalised within $\RCA_0$, and used to prove that
$\mathrm{HB}$ implies weak König's lemma \citep[\S IV.1]{Simpson2009}.


\section{Finitism and conservativity}
\label{sec:finitism}

As far as purely number-theoretic statements are concerned, $\WKL_0$ and
$\RCA_0$ prove exactly the same ones: if $\varphi$ is an arithmetical sentence
such that $\WKL_0 \vdash \varphi$, then $\RCA_0 \vdash \varphi$
\citep[\S IX.2]{Simpson2009}.
Moreover, by a theorem of Parsons, $\WKL_0$ and $\RCA_0$ are both conservative
for $\Pi^0_2$ sentences over the theory $\PRA$ of primitive recursive
arithmetic \citep{Parsons1970,Friedman1976}.
These facts led \citet{Simpson1985, Simpson1988} to suggest that the mathematics
provable in $\WKL_0$, and other systems which are conservative over $\PRA$ for
$\Pi^0_1$ sentences, is \emph{finitistically reducible}, and that this
constitutes a partial realisation of Hilbert's programme.

To understand this claim we first need to briefly reexamine Hilbert's programme.
Hilbert took finitary mathematics---the part of mathematics with finitary
content---to include number-theoretic equations and inequalities, but also
universal generalisations of such formulas, i.e.\ $\Pi^0_1$ sentences.
Hilbert aimed to justify infinitary mathematics in a finitary way. Initially the
intent of his programme was to provide finitary proofs of the consistency of
infinitary theories. In the second half of the 1920s its aim shifted somewhat,
towards proving the conservativity of infinitary theories over finitary
theories. This new version of the programme can be understood in terms of an
analogy with physics suggested by \citet{Weyl1925}. The finitary or `real' part
of a mathematical theory is accorded a status analogous to that of the
observation sentences of a physical theory, the part of the theory that describe
the outcomes of empirical observations. If the infinitary or `ideal' part of the
theory, analogous to the theoretical sentences of a physical theory, imply a
finitary statement, then that finitary statement must already be provable using
finitary means alone, just as observation sentences implied by a physical theory
must be (in principle) verifiable by experiment. In other words, a mathematical
theory must be conservative over its finitary part.

In an influential article on Hilbertian finitism, \citet{Tait1981} argues for
two theses. Firstly, the finitist functions $f : \N^k \to \N$ are precisely the
primitive recursive functions. Secondly, the finistically provable part of
mathematics consists of those $\Pi^0_1$ sentences which are provable in the
formal system $\PRA$ of primitive recursive arithmetic.%
\footnote{%
Tait's theses have gained broad but not universal acceptance.
\citet{Kreisel1958,Kreisel1970} considers them too restrictive, and argues that
the finistically provable sentences include those provable in the stronger
system $\PA$ of first-order Peano arithmetic.
\citet{Ganea2010} instead argues that Tait's thesis is too permissive, and that
finitistic provability includes no more than sentences provable from the
equational theory of Kalmár elementary functions.
}
Building on Tait's theses, Simpson claims that the $\Pi^0_1$ conservativity of
$\RCA_0$ and $\WKL_0$ over $\PRA$ shows that Hilbert's programme can succeed, at
least in part, because the infinitary mathematics provable in these systems is
so substantial. This argument appears to be bolstered by the fact that the
conservativity of these systems over $\PRA$ is finitistically provable, in the
following sense. \citet{Sieg1985a} constructed a primitive recursive function
$g : \N \to \N$ which, given as input a proof $p$ of a $\Pi^0_2$ sentence
$\varphi$ in $\WKL_0$, produces as output a proof $g(p)$ of $\varphi$ in $\PRA$.
Because the conservativity theorem is, per Tait's thesis, finitistically
provable, Simpson argues that Hilbert's programme is successfully realised for
the infinitary mathematics provable in $\WKL_0$.

As \citet{Burgess2010} points out, a closer reading of the situation suggests a
difficulty with Simpson's argument. Suppose that $\varphi$ is a $\Pi^0_1$
sentence, hence finitistically meaningful, with a proof $p$ in $\WKL_0$. Then
$g(p)$ is a finitist proof per Tait's thesis: all the finitist has to do is
check that all of the axioms used in the proof are finitistically acceptable,
which they are, since they are axioms of $\PRA$. The problem is that this is
not sufficient for the finitist to have confidence in the infinitary theory
$\WKL_0$ as a whole. For that, they would need not only to know (in a
finitistically provable way, i.e.\ in $\PRA$) that $\WKL_0$ is $\Pi^0_1$
conservative over $\PRA$, but that any $\Pi^0_1$ consequence of $\WKL_0$ is
finitistically true. In other words, they would need to know the $\Pi^0_1$
reflection principle for $\WKL_0$, the scheme
\begin{equation*}
    \mathrm{Prov}_{\WKL_0}(\gcode{\varphi}) \rightarrow \varphi
\end{equation*}
for all sentences $\varphi \in \Pi^0_1$.
The conservativity theorem only guarantees, for the finitist, that
\begin{equation*}
    \mathrm{Prov}_{\WKL_0}(\gcode{\varphi})
    \rightarrow
    \mathrm{Prov}_{\PRA}(\gcode{\varphi}),
\end{equation*}
since to establish $\varphi$ as finitistically provable we need to invoke Tait's
thesis, the claim that
\begin{equation*}
    \mathrm{Prov}_{\PRA}(\gcode{\varphi}) \rightarrow \varphi
\end{equation*}
for any $\Pi^0_1$ sentence $\varphi$.
This scheme is not provable in $\PRA$, since an instance of the scheme is
\begin{equation*}
    \mathrm{Prov}_{\PRA}(\gcode{0=1}) \rightarrow 0=1,
\end{equation*}
which is equivalent to the consistency statement for $\PRA$, and hence
unprovable in $\PRA$ by Gödel's second incompleteness theorem.%
\footnote{%
See \citet{Giaquinto1983} and \citet{Dean2015} for related discussion.
}

The problem is that Tait's thesis is a claim made from outside the finitary
standpoint, not from within it \citep[p.~527]{Tait1981}. Finitists can know of
individual infinitary proofs that they establish finitistically true statements,
via a finitarily provable conservativity theorem, since such conservativity
theorems deliver a finite proof object that can be verified by the finitist to
only use finitistically acceptable axioms. But they cannot know, at least for
systems such as $\RCA_0$ or $\WKL_0$, that these systems only prove
finitistically true $\Pi^0_1$ sentences, since this amounts to knowing the
consistency of $\PRA$, something which is unavailable to them.


\section{Predicativity and definability}
\label{sec:predicativism}

In the wake of the class-theoretic paradoxes, Russell realised that
propositional functions do not always define a corresponding class. He called
\emph{predicative} those propositional functions $\varphi(x)$ such that the
corresponding class $\left\{ x : \varphi(x) \right\}$ exists. Propositional
functions like $x \not\in x$ which do not, on pain of contradiction, define a
class, Russell called \emph{impredicative} \citep{Russell1907a}. This
distinction gave way to one framed in terms of the theory of types, and of the
`vicious circle principle': ``no totality can contain members defined in terms
of itself'' \citep[p.~237]{Russell1908}. On this account, a propositional
function $\varphi(x)$ has type $n+1$ just in case its arguments and values have
type $n$, and all variables bound by quantifiers in $\varphi(x)$ have type $n$
or lower. A predicative function is one which appears in this hierarchy, i.e.\
has type $n+1$ for some $n$, where $0$ is the type of individuals.%
\footnote{%
The history of Russell's thought is traced in \citet{IrvineDeutsch2020}.
}

\citet{Poincare1908} objected to the details of Russell's attempts to spell out
a solution, but agreed that the source of the paradoxes lay in definitions which
contained vicious circles, in the sense that if one defines two concepts $C$ and
$C'$, their definitions are impredicative if the concept $C$ occurs in the
definition of $C'$, and conversely. However, Poincaré also calls definitions
predicative if they are ``not changed by the introduction of new elements''
\citep{Poincare1910}.
\citet{Weyl1918} drew on the ideas of Russell and Poincaré in developing his
own view, one in which vicious circles were to be avoided by restricting
comprehension to predicates which quantified only over natural numbers. Despite
the fact that Weyl swiftly moved on from predicativism, adopting Brouwer's
intuitionism, his approach exerted a significant influence on later developments
in predicative analysis and on reverse mathematics. Dedekind had already shown
that fundamental theorems of 19th century analysis were equivalent to one
another and to the least upper bound principle \citep[\S VII]{Dedekind1872}.
Weyl now showed that versions of these theorems could be proved in a predicative
framework.

After a fallow period, predicativity underwent a revival in the mid-1950s,
becoming connected to the emerging subject of definability theory.
\citet{Grzegorczyk1955,Kondo1958}, and \citet{Mostowski1959} all suggested that
predicative analysis in the sense of Weyl could be identified with elementary
analysis, in which every set of natural numbers is given by an elementary
definition, i.e.\ one in the first-order language of arithmetic. This amounts to
studying analysis in the $\omega$-model
$\ARITH = \left\{
    X \subseteq \omega : X \text{ is arithmetically definable}
\right\}$.
Grzegorczyk and Mostowski both provided axiomatisations of elementary analysis,
which can be seen as ancestors of the theory $\ACA_0$. The axioms of this
system can all be justified on predicative grounds, since it postulates the
existence only of sets which are definable in terms of objects that are already
`given', namely the natural numbers. Theorems provable in $\ACA_0$ are thus
predicatively provable.

The fundamental idea of predicative definability, in which sets can be defined
in terms of previously given or already existing objects, is prima facie more
general than arithmetical definability, since it seems predicatively acceptable
to then iterate the procedure, taking the arithmetical sets as given and
defining new sets on that basis.
One way of doing this is by iterating arithmetical definability into the
transfinite, along well-orderings which are predicatively acceptable, and
\citet{Wang1954} and \citet{Lorenzen1955} made early forays into this approach.
This idea underlaid the proposal of \citet{Kreisel1960} to identify the
predicatively definable sets of natural numbers with those occurring in the
hyperarithmetical hierarchy.%
\footnote{%
The hyperarithmetical sets are those which are computable relative to the
$\alpha$th iteration of the Turing jump operator starting from the empty set,
where $\alpha$ is a computable ordinal \citep[Part~A]{Sacks1990}.
}
On this basis, Kreisel proposed several axiom systems true in the $\omega$-model
$\HYP = \left\{ X \subseteq \omega : \text{$X$ is hyperarithmetical} \right\}$,
including $\DiCA_0$ and $\SinAC{1}_0$.
These systems of hyperarithmetical analysis did not prove substantially more
mathematically fruitful than elementary analysis, with Kreisel writing that
``in the portions of analysis developed by working mathematicians, a theorem is
either derivable by means of the arithmetic comprehension axiom or else not
predicative at all'' \citep[p.~316]{Kreisel1962}.
Hyperarithmetical analysis has only recently started to yield theorems which are
not arithmetically true, such as statements concerning indecomposable linear
orders \citep{Montalban2006}, and versions of Halin's infinite ray theorem
\citep*{BarnesGohShore2022}.

The two remaining Big Five systems, $\ATR_0$ and $\PiCA_0$, both have a level of
stability and mathematical fruitfulness not achieved by theories of
hyperarithmetical analysis. However, on the widely accepted thesis that the
ordinal $\Gamma_0$ marks the outer limit of predicativity, both of these
theories are impredicative. This is most obvious in the case of $\PiCA_0$, since
its proof-theoretic ordinal is substantially larger than $\Gamma_0$, and its
characteristic axiom is prima facie impredicative.%
\footnote{%
This analysis of the limits of predicativity is surveyed by
\citet{Feferman2005}.
The proof-theoretic ordinal of a system $S$ is usually understood as the least
ordinal $\alpha$ such that $S$ cannot prove that $\alpha$ is well-ordered,
modulo a reasonable computable presentation (an ordinal notation system) of
$\alpha$ and ordinals $\beta < \alpha$.
Systems such as $\ATR_0$ which are not predicative themselves, but have
proof-theoretic ordinal $\Gamma_0$, are sometimes called \emph{predicatively
reducible}. \citet{Simpson1985} argues that $\ATR_0$ and other predicatively
reducible systems can be used instrumentally by the predicativist. However, this
analysis seems vulnerable to an argument parallel to that deployed against
finististic reducibility, as suggested by \citet[p.~140]{Burgess2010} and
discussed by \citep[\S 5.5]{Eastaugh2024}.
The $\Gamma_0$ analysis of the limits of predicative provability is not
universally accepted, with \citet{Weaver2009} a notable dissenter.
}
\citet{Kreisel1959} constructed a recursive counterexample to the
Cantor--Bendixson theorem, which states that every closed set is the union of a
perfect set and a countable set. Kreisel's construction is of a computable (code
for a) closed set $C$ such that, by the Cantor--Bendixson theorem,
$C = P \cup S$ with $P$ a perfect set and $S$ a countable set, but such that
both $P$ and $S$ are coanalytic ($\Pi^1_1$) but not analytic ($\Sigma^1_1$). It
follows that neither are hyperarithmetical, and hence not predicative by
Kreisel's analysis of predicative definability. 
Kreisel's argument can be internalised within $\ACA_0$ to show that the
Cantor--Bendixson theorem implies $\Pi^1_1$ comprehension
\citep[\S VI.1]{Simpson2009}. The abundance of reversals to $\PiCA_0$ and
$\ATR_0$ suggests that these systems form natural accumulation points in the
hierarchy of subsystems of second-order arithmetic, even in the absence of
justifications for these systems in terms of historical foundational programmes.


\section{The intrinsic significance of reversals}
\label{sec:intrinsic}

Philosophy of mathematics is broader than foundations of mathematics, especially
when one considers this in the narrow sense of the early 20th century
foundational programmes. One might reasonably wonder what reverse mathematics
can contribute to philosophy of mathematics more broadly, beyond measuring the
limits of the mathematics that can be developed from within a given foundational
perspective, especially from a more broadly realist point of view on which the
axioms of subsystems of second-order arithmetic are viewed as true, but
inadequate to formalise all of mathematics. We can start to understand what is
at stake here by asking: What do we learn from reversals?

The received view, as presented most influentially by
\citet[pp.~1--2]{Simpson2009}, is that reversals show which \emph{set existence
principles} are needed to prove the theorem in question. The use of this term
emphasises issues of mathematical ontology, and might lead one to infer that the
practice of reverse mathematics involves some underlying scepticism about the
reliability of set-theoretic methods, or the ontological picture of the
cumulative hierarchy of sets which is conventionally used to justify them.
This is far from being the case: most work in reverse mathematics is neither
motivated by a particular foundational programme, nor constrained by a desire to
do away with tools and concepts that go beyond those formalisable in
second-order arithmetic. It does, however, leave open how best to answer the
question of what we learn from reversals, as a general issue in the philosophy
of mathematics.

\citet{Hirschfeldt2014} has argued that reversals reveal the \emph{combinatorial
core} of a theorem, the combinatorial principle that is essential to or
underpins any proof of that theorem, regardless of whether it makes an explicit
appearance in that proof. For example, the combinatorial core of Lindenbaum's
lemma is weak König's lemma \citep[p.~9]{Hirschfeldt2014}.
The use of $\mathrm{WKL}$ is clear in the standard proof of Lindenbaum's lemma.
Fix an enumeration $\str{\varphi_n : n \in \N}$ of the sentences of a countable
language $\lang$, and for any $\varphi_i$ in the enumeration let $\varphi_i^1
\equiv \varphi$ and $\varphi_i^0 \equiv \neg\varphi$. Given an $\lang$-theory
$S$ and a finite sequence $\sigma \in 2^{<\N}$, set
\begin{equation*}
    S_\sigma = S \cup \left\{
        \varphi_i^{\sigma(i)} : i < \lh{\sigma}
    \right\}.
\end{equation*}
We form a binary tree $T_S$ by putting $\sigma \in T_S$ if there is no proof of a
contradiction of length $\leq \lh{\sigma}$ from $S_\sigma$. This is a computable
process, so $T_S$ can be proved to exist in $\RCA_0$, relative to $S$. $T_S$ is
infinite if and only if $S$ is consistent, so by $\mathrm{WKL}$ there is an
infinite path $P : \N \to 2$ through $T_S$. The set
$S^* = \left\{ \varphi_i : P(i) = 1 \right\}$
is a maximal consistent set extending $S$. Weak König's lemma is essential
here, since there are computable consistent sets $T$ (such as $\PA$) which
cannot be extended to a computable maximal consistent set $S^*$. This points to
a crucial aspect of reverse mathematics, namely the fact that the characteristic
axioms of the Big Five have a computability-theoretic character.

This is most easily grasped by examining the $\omega$-models of the Big Five. We
have already discussed some examples: the computable sets $\REC$ form an
$\omega$-model of $\RCA_0$, the arithmetical sets $\ARITH$ form an
$\omega$-model of $\ACA_0$, and so on. These examples are quite representative,
and reveal important structural properties common between all $\omega$-models of
those systems. For example, every $\omega$-model $\mathcal{X} \subseteq
\powset{\omega}$ of $\RCA_0$ is downwards closed under Turing reducibility,
meaning that if $X \in \mathcal{X}$ and $Y$ is computable by a Turing machine
with an oracle for $X$, then $Y \in \mathcal{X}$. Similarly every $\omega$-model
$\mathcal{X} \subseteq \powset{\omega}$ of $\ACA_0$ is closed under arithmetical
reducibility: if $X \in \mathcal{X}$ and $Y$ is computable from $X$ via a finite
number of Turing jumps, then $Y \in \mathcal{X}$.%

\citet{Eastaugh2019}, following other authors such as \citet{Kreisel1968,
Shore2010}, and \citet*{ChongFengSlaman2014}, foregrounds the computability-%
theoretic nature of the characteristic axioms of the Big Five, and argues that
they express \emph{closure conditions} on the second-order part of the model.
This suggests the following unified way of looking at what we learn from
reversals. Stronger axioms express stronger closure conditions, which guarantee
that more sets exist (the \emph{existential} aspect), that stronger
combinatorial properties hold (the \emph{combinatorial} aspect), and that models
of the stronger axioms are closed under stronger computability-theoretic
reducibility relations (the \emph{computational} aspect). They are thus
intertwined facets of the same basic phenomenon.

Reversals thus play an explanatory role, revealing the logico-combinatorial-%
computational \emph{structure} that supports (and is supported by) a particular
theorem---or, considered slightly differently, as revealing the logico-%
combinatorial-computational \emph{content} of that theorem.%
\footnote{%
Although this explanatory role can be partially understood in semantic terms, it
seems implausible that reversals can give a satisfying theory of mathematical
content broadly understood. Such a view might identify the content of a theorem
$\varphi$ with its set of models, or perhaps with the equivalence class of
statements which are provably equivalent to $\varphi$ in $\RCA_0$. But as we
have seen, these equivalence classes can be very large, and it seems implausible
to say that for the practising analyst the Bolzano--Weierstra\ss{} theorem has
the same content as König's lemma. See \citet{AranaMancosu2012} for some
related considerations in the context of planar and solid geometry.
}
For example, a theory adequate to the analysis of Weierstra\ss{} and Dedekind
must include the Bolzano--Weierstra\ss{} theorem and the Cauchy convergence
theorem. The fact that we can reason backwards from those theorems to the axiom
scheme of arithmetical comprehension, as well as forwards from $\ACA_0$ to those
theorems, is revealing of both an important logical and combinatorial structure
in the theory of analysis, and the computability-theoretic content of those
theorems. This latter fact was not evident from the practice of real analysis
alone, and needed to be shown by work in mathematical logic.
It is equally revealing that although the sequential Heine--Borel theorem is
derivable in this system, one cannot reverse this implication, and the models
of $\WKL_0$ (the weaker system necessary and sufficient to prove the sequential
Heine--Borel theorem) have quite different properties to those of $\ACA_0$.
Compactness principles like Heine--Borel thus belong to a quite different
species than completeness principles like the monotone convergence or
Bolzano--Weierstra\ss{} theorems. This is so even though the adoption of
completeness principles was historically key to proving the Heine--Borel
theorem.%
\footnote{%
See \citet{AndreEngdahlParker2013} and the extensive discussion by
\citealt[pp.~20--25]{Hallett1979}.
}


\section{Reverse mathematics and nominalism}
\label{sec:nominalism}

The indispensability argument seeks to justify realism about mathematics on the
basis of mathematics' role in our best scientific theories, where those theories
are construed as true or at least truth-apt claims about the nature of the
world. Mathematical entities, on Quine's view, are epistemically on a par with
other theoretical entities held to exist by scientific theories, since they are
indispensable to our best science.
A natural question which emerges in this context is how much mathematics is
indispensable to science. Quine took even the irrational numbers to stand in
need of justification, in terms of the ways in which admitting them into our
theories simplifies computations and generalisations---in other words,
theoretical virtues that go beyond the mere prediction of observed values in
possible experiments. Set theory beyond the real numbers is admissible only
because it completes and systematises our best theories of applied mathematics
\citep[p.~400]{Quine1998}.

Putnam, another exponent of the indispensability argument, held that physics
could survive---if not thrive---using only the mathematical resources of
predicative set theory. On the other hand, the mathematical usefulness of
impredicative set theory for the mathematics ultimately applied in physics
constitutes an argument for its truth, albeit not as strong as that for the
indispensable predicative theory \citep[p.~55--56]{Putnam1971}. The first part
of Putnam's view has been embraced by those with nominalist or constructivist
sympathies, since if the indispensability argument does not justify all or even
much of classical set theory, then the door might be open to reconstructing the
scientifically applicable part of mathematics in a nominalist or broadly
constructive way. It offers an obvious role for reverse mathematics, in
characterising the mathematical axioms necessary for deriving theorems which
are indispensably applied in physics, or biology, or even in social sciences
such as economics.

\citet{Feferman1988,Feferman1992} has defended a view along these lines, holding
that the mathematics that is indispensable to our best science is predicative.
The body of results proved in reverse mathematics endows Feferman's claim with a
certain plausibility, because a substantial fragment of 19th and 20th century
analysis can be recovered in predicative systems such as $\ACA_0$. This includes
many theorems which are widely used in scientific applications.
A counterexample to this claim would have to have the form of a mathematical
theorem which is indispensable to some scientific application but unprovable in
$\ACA_0$ (or some conservative extension thereof). One suggestion has been that
non-separable spaces are used in nontrivial ways in applied mathematics,
especially physics. Feferman's response \citep[p.~281]{Feferman1998} is that
although such spaces are widely used, it is not clear that the non-separability
of the spaces really makes a difference to the applications in question. This
remains an area where careful, detailed work in logic and the foundations of
physics remains necessary.

More strictly nominalist authors including \citet{Hellman1989,Hellman1999} and
\citet{Bueno2001} have sought to use reverse mathematical results to defuse the
indispensability argument. Their general strategy is to use reversals to
determine the scope of scientifically applicable mathematics in terms of
subsystems of second-order arithmetic, and then argue that those subsystems can
be reinterpreted in a nominalistically acceptable fashion.
\citet{Hellman1999}, for example, argues that one can develop a nominalistic
account of the natural numbers, and thus an indispensability argument is not
needed in order to justify predicative theories such as $\ACA_0$. This is
because for such theories, membership ascriptions of the form $n \in X$ can be
reinterpreted as ascriptions of the form $\varphi(n)$ where $\varphi$ is an
arithmetical formula defining the set $X$. These arithmetical ascriptions are
taken to be ontologically non-committing due to the claimed availability of a
nominalisation of the natural numbers. Where quantification over sets is needed,
it can be replaced by quantification over arithmetical formulas. Note however
that the acceptability of the nominalisation strategy rests on the complexity of
the sets which the subsystem of second-order arithmetic proves to exist: a
strategy that works for theorems provable in $\ACA_0$ will not directly transfer
to theorems provable in $\PiCA_0$. Arithmetical transfinite recursion and
stronger axioms like $\Pi^1_1$ comprehension are therefore the ``next natural
target of indispensability arguments'' \citep[p.~37]{Hellman1999}.


\section{The indispensability of mathematics to philosophy}
\label{sec:reverse_phil}

We close by briefly considering a different way in which mathematics may be
indispensable, namely to our philosophical and conceptual work, rather than to
empirical science. A self-contained example is Kreisel's
\citeyearpar{Kreisel1967a} \emph{squeezing argument} that the informal concept
of logical validity for first-order sentences, $\mathrm{Val}(x)$, is
co-extensional with the formal concept of logical validity for first-order
sentences, i.e.\ the set of first-order sentences $\varphi$ which are true in
all set-sized models $\mathcal{M}$. Kreisel aims to establish that the informal
concept is determinate, with a sharp boundary that excludes the possibility of a
first-order sentence whose validity is unsettled. He does this by arguing that
$\mathrm{Val}(x)$ has the following properties.
\begin{enumerate}
    \item[(P1)] For all $\varphi$, if $\varphi$ is derivable in classical
        first-order logic, then $\mathrm{Val}(\varphi)$.
    \item[(P2)] For all $\varphi$, if $\mathrm{Val}(\varphi)$,
        then $\mathcal{M} \models \varphi$ for all set-theoretic models
        $\mathcal{M}$.
\end{enumerate}
The third premise for the argument is Gödel's completeness theorem for
first-order logic.
\begin{enumerate}
    \item[(P3)] For all $\varphi$, if $\mathcal{M} \models \varphi$ for all
        set-theoretic models $\mathcal{M}$, then $\varphi$ is derivable in
        classical first-order logic.
\end{enumerate}
This logically entails the conclusion,
\begin{enumerate}
    \item[(C)] For all $\varphi$, $\mathrm{Val}(\varphi)$ if and only if
        $\mathcal{M} \models \varphi$ for all set-theoretic models
        $\mathcal{M}$.
\end{enumerate}
In other words, informal validity coincides with formal validity.%
\footnote{%
A similar argument appears, without citing Kreisel, in Quine's discussion of
logical truth in \citep[ch.~4]{Quine1986}.
\citet{Field1989} appears to have coined the term `squeezing argument'.
}

Two of the argument's premises are philosophical or conceptual in nature,
concerning the content of our informal conception of validity. The third appears
purely mathematical, but it is nevertheless essential to the soundness of the
argument: if P3 is false then the conclusion is not necessitated by the other
premises.
One way of understanding this is as follows. $\mathrm{Con}(\PA)$ is a sentence
of first-order number theory, so assuming that $\PA$ is in fact consistent,
$\mathrm{Con}(\PA)$ is true in the standard model of arithmetic. As an
$\omega$-model, $\REC \models \mathrm{Con}(\PA)$, but there is no $X \in \REC$
which codes a countable model of $\PA$, since all such sets are non-computable.
In other words, Gödel's completeness theorem is not computably true, and in
fact it is equivalent over $\RCA_0$ to weak König's lemma
\citep[\S IV.3]{Simpson2009}.
From this we can conclude that in accepting Kreisel's argument---that is, in
accepting its premises---we must also accept the axioms of $\WKL_0$ and the
existence of non-computable sets.%
\footnote{%
A form of this idea is implicit in \citet[p.~148]{Etchemendy1990}.
}

\citet{DeanKurokawa2026} consider this in the context of a larger project of
analysing the methodology of what Kreisel called `informal rigour'. This seems
to offer the prospect of discovering further examples of the indispensability of
mathematics to philosophy, by examining the mathematical premises of other
squeezing arguments.%
\footnote{%
There are many antecedents to this general approach, in which a mathematical
theorem is used to close a circle of containments involving different informal
mathematical concepts. It concerns a form of argument which one might consider
native to philosophy of mathematics, in which the fact that a given theorem
which bears on a philosophical analysis is provable (or not provable) in a given
deductive system.
Examples include Skolem's reworking of the downward Löwenheim--Skolem theorem
\citep{Skolem1920} as the countable model theorem \citep{Skolem1923a} in order
to eliminate the use of the axiom of choice, or the status of the completeness
theorem for first order logic in programmes such as Field's
\citep{Field2016,Field1991}.
The specific idea of an indispensability argument for mathematical realism based
on the use of mathematics in philosophy emerged in discussions between Marianna
Antonutti Marfori, Walter Dean, and the author.
}
However, squeezing arguments are not the only cases in which we can apply the
general methodology of using reverse mathematics to analyse the strength of the mathematical premises of philosophical arguments. A recent example is Dean and
Sanders's \citeyearpar{DeanSanders2025} reverse mathematical analysis of the
continuous sorites paradox introduced by \citet{WeberColyvan2010}. They show
that an $\lang_2$ formalisation of the mathematical principles necessary to
generate the paradox reverses to arithmetical comprehension. One way of
understanding their result is that the computable real numbers can contain
continuous sorites sequences without paradox, while another is that the
unknowability of the location of a sharp boundary of where such sequences end
could correspond to the non-computability of such a boundary, suggesting the
possibility of a computational variant on epistemicism. Studying this
mathematically entangled paradox through the lens of reverse mathematics thus
reveals new perspectives on the problem, including introducing a novel sorites
argument mediated by the Heine--Borel theorem \citep[\S 5]{DeanSanders2025}.

A final example is Arrow's theorem, a result in voting theory which is sometimes
taken to show (e.g.\ by \citealt{Riker1982}) that populist democracy is
impossible. \citet{Eastaugh2024a} shows that Arrow's theorem is provable in
primitive recursive arithmetic, and thus will be acceptable to not just
platonists but a wide variety of nominalist and constructive perspectives. On
the other hand, Arrow's theorem is known to fail when the set of voters is
infinite, by Fishburn's \citeyearpar{Fishburn1970} possibility theorem. Given
the contemporary interest in infinite societies in population ethics, modelling
societies which can continue indefinitely into the future, Fishburn's theorem
might be taken to suggest that infinite societies differ from finite ones in
admitting democratic choice procedures. However, Fishburn's theorem is
equivalent over $\RCA_0$ to arithmetical comprehension
\citep[theorem~5.4]{Eastaugh2024a}. It follows that there are computable
societies with no computable non-dictatorial social welfare function. Although
such social welfare functions exist mathematically, they therefore function as
democratic decision-making procedures only in an idealised sense, even in the
case of infinite societies in which individuals' preferences remain computable.

This essay has aimed to present reverse mathematics as an intrinsically
interesting research programme in contemporary mathematical logic. In addition
to its connection to traditional foundational programmes in the philosophy of
mathematics, it also provides a means of refining contemporary debates about
mathematical content, realism, and indispensability. More recently it has
entered the discussion of purity of methods \citep{Arana2024} and the relation
between consistency and existence \citep{Dean2021a}. The future therefore looks
bright for further philosophical engagement with the methods and results of
reverse mathematics.


\bibliographystyle{abbrvnat}
\bibliography{reverse-math-blackwell}

\begin{thebibliography}{103}
\providecommand{\natexlab}[1]{#1}
\providecommand{\url}[1]{\texttt{#1}}
\expandafter\ifx\csname urlstyle\endcsname\relax
  \providecommand{\doi}[1]{doi: #1}\else
  \providecommand{\doi}{doi: \begingroup \urlstyle{rm}\Url}\fi

\bibitem[Andre et~al.(2013)Andre, Engdahl, and Parker]{AndreEngdahlParker2013}
N.~R. Andre, S.~M. Engdahl, and A.~E. Parker.
\newblock An analysis of the first proofs of the {Heine}--{Borel} theorem.
\newblock \emph{Loci}, 4, August 2013.
\newblock \doi{10.4169/loci003890}.

\bibitem[Arana(2024)]{Arana2024}
A.~Arana.
\newblock \emph{Elements of Purity}.
\newblock Elements in the Philosophy of Mathematics. Cambridge University
  Press, Cambridge, 2024.
\newblock \doi{10.1017/9781009052719}.

\bibitem[Arana and Mancosu(2012)]{AranaMancosu2012}
A.~Arana and P.~Mancosu.
\newblock On the relationship between plane and solid geometry.
\newblock \emph{The Review of Symbolic Logic}, 5\penalty0 (2):\penalty0
  294--353, June 2012.
\newblock \doi{10.1017/S1755020312000020}.

\bibitem[Barnes et~al.(2022)Barnes, Goh, and Shore]{BarnesGohShore2022}
J.~S. Barnes, J.~L. Goh, and R.~A. Shore.
\newblock Theorems of hyperarithmetic analysis and almost theorems of
  hyperarithmetic analysis.
\newblock \emph{The Bulletin of Symbolic Logic}, 28\penalty0 (1):\penalty0
  133--149, March 2022.
\newblock \doi{10.1017/bsl.2021.70}.

\bibitem[Bishop(1967)]{Bishop1967}
E.~Bishop.
\newblock \emph{Foundations of Constructive Analysis}.
\newblock McGraw-Hill, New York, 1967.

\bibitem[Bishop(1975)]{Bishop1975}
E.~Bishop.
\newblock The crisis in contemporary mathematics.
\newblock \emph{Historia Mathematica}, 2\penalty0 (4):\penalty0 507--517,
  November 1975.
\newblock \doi{10.1016/0315-0860(75)90113-5}.

\bibitem[Bueno(2001)]{Bueno2001}
O.~Bueno.
\newblock Logicism revisited.
\newblock \emph{Principia}, 5\penalty0 (1--2):\penalty0 99--124, 2001.
\newblock URL \url{https://dialnet.unirioja.es/servlet/articulo?
  codigo=5251158}.

\bibitem[Burgess(2010)]{Burgess2010}
J.~P. Burgess.
\newblock {On the outside looking in: a caution about conservativeness}.
\newblock In S.~Feferman, C.~Parsons, and S.~G. Simpson, editors, \emph{Kurt
  G{\"o}del, Essays for His Centennial}, number~33 in Lecture Notes in Logic,
  pages 128--141. Cambridge University Press, 2010.
\newblock \doi{10.1017/CBO9780511750762.009}.

\bibitem[Cantor(1872)]{Cantor1872}
G.~Cantor.
\newblock {Ueber die Ausdehnung eines Satzes aus der Theorie der
  trigonometrischen Reihen}.
\newblock \emph{Mathematische Annalen}, 5\penalty0 (1):\penalty0 123--132,
  1872.
\newblock \doi{10.1007/BF01446327}.

\bibitem[Chong et~al.(2014)Chong, Feng, Slaman, and
  Woodin]{ChongFengSlaman2014}
C.~Chong, Q.~Feng, T.~A. Slaman, and W.~H. Woodin, editors.
\newblock \emph{Infinity and Truth}, volume~25 of \emph{Lecture Notes Series,
  Institute for Mathematical Sciences, National University of Singapore}.
\newblock World Scientific, Singapore, 2014.
\newblock \doi{10.1142/9789814571043_fmatter}.

\bibitem[Dauben(1979)]{Dauben1979}
J.~W. Dauben.
\newblock \emph{Georg Cantor: His Mathematics and Philosophy of the Infinite}.
\newblock Princeton University Press, Princeton, New Jersey, 1979.

\bibitem[Dean(2015)]{Dean2015}
W.~Dean.
\newblock Arithmetical reflection and the provability of soundness.
\newblock \emph{Philosophia Mathematica}, 23\penalty0 (1):\penalty0 31--64,
  February 2015.
\newblock \doi{10.1093/philmat/nku026}.

\bibitem[Dean(2020)]{Dean2021a}
W.~Dean.
\newblock On consistency and existence in mathematics.
\newblock \emph{Proceedings of the Aristotelian Society}, 120\penalty0
  (3):\penalty0 349--393, October 2020.
\newblock \doi{10.1093/arisoc/aoaa017}.

\bibitem[Dean and Kurokawa(2026)]{DeanKurokawa2026}
W.~Dean and H.~Kurokawa.
\newblock On the methodology of informal rigour: set theory, semantics, and
  intuitionism.
\newblock \emph{Journal of Philosophical Logic}, 2026.
\newblock \doi{10.1007/s10992-026-09839-5}.

\bibitem[Dean and Sanders(2025)]{DeanSanders2025}
W.~Dean and S.~Sanders.
\newblock From real analysis to the continuous sorites via reverse mathematics.
\newblock \emph{The Review of Symbolic Logic}, 18\penalty0 (3), 2025.
\newblock \doi{10.1017/S1755020325000061}.

\bibitem[Dean and Walsh(2017)]{DeanWalsh2017}
W.~Dean and S.~Walsh.
\newblock The prehistory of the subsystems of second-order arithmetic.
\newblock \emph{The Review of Symbolic Logic}, 10\penalty0 (2):\penalty0
  357--396, 2017.
\newblock \doi{10.1017/S1755020316000411}.

\bibitem[Dedekind(1872)]{Dedekind1872}
R.~Dedekind.
\newblock \emph{Stetigkeit und irrationale Zahlen}.
\newblock Vieweg, 1872.
\newblock English translation in {\citet{Dedekind1901}}.

\bibitem[Dedekind(1888)]{Dedekind1888}
R.~Dedekind.
\newblock \emph{Was sind und was sollen die Zahlen?}
\newblock Vieweg, 1888.
\newblock English translation in {\citet{Dedekind1901}}.

\bibitem[Dedekind(1901)]{Dedekind1901}
R.~Dedekind.
\newblock \emph{Essays on the Theory of Numbers}.
\newblock Open Court, 1901.
\newblock English translations of {\citet{Dedekind1872}} and
  {\citet{Dedekind1888}}, edited and translated by W. W. Beman.

\bibitem[{Dzhafarov} and Mummert(2022)]{DzhafarovMummert2022}
D.~D. {Dzhafarov} and C.~Mummert.
\newblock \emph{Reverse Mathematics: Problems, Reductions, and Proofs}.
\newblock Number~7 in Theory and Applications of Computability. Springer, Cham,
  2022.
\newblock \doi{10.1007/978-3-031-11367-3}.

\bibitem[Eastaugh(2019)]{Eastaugh2019}
B.~Eastaugh.
\newblock Set existence principles and closure conditions: unravelling the
  standard view of reverse mathematics.
\newblock \emph{Philosophia Mathematica}, 27\penalty0 (2):\penalty0 153--176,
  June 2019.
\newblock \doi{10.1093/philmat/nky010}.

\bibitem[Eastaugh(2024)]{Eastaugh2024}
B.~Eastaugh.
\newblock Reverse mathematics.
\newblock In E.~N. Zalta and U.~Nodelman, editors, \emph{The {Stanford}
  Encyclopedia of Philosophy}. Metaphysics Research Lab, Stanford University,
  summer 2024 edition, 2024.
\newblock URL \url{https://plato.stanford.edu/archives/sum2024/entries/reverse-
  mathematics/}.

\bibitem[Eastaugh(2025)]{Eastaugh2024a}
B.~Eastaugh.
\newblock {Arrow}'s theorem, ultrafilters, and reverse mathematics.
\newblock \emph{The Review of Symbolic Logic}, 18\penalty0 (2):\penalty0
  439--462, June 2025.
\newblock \doi{10.1017/S1755020324000054}.

\bibitem[Etchemendy(1990)]{Etchemendy1990}
J.~Etchemendy.
\newblock \emph{The Concept of Logical Consequence}.
\newblock Harvard University Press, 1990.

\bibitem[Ewald(1996)]{Ewald1996b}
G.~Ewald.
\newblock \emph{Combinatorial Convexity and Algebraic Geometry}.
\newblock Springer, New York, 1996.
\newblock \doi{10.1007/978-1-4612-4044-0}.

\bibitem[Feferman(1988)]{Feferman1988}
S.~Feferman.
\newblock {Weyl} vindicated: {\emph{Das Kontinuum}} seventy years later.
\newblock In \emph{Termi e prospettive della logica e della filosofia della
  scienza contemporanee, vol. 1}, pages 59--93. CLUEB, 1988.
\newblock Reprinted with minor corrections and additions in
  {\citet[pp.~249--283]{Feferman1998}}.

\bibitem[Feferman(1992)]{Feferman1992}
S.~Feferman.
\newblock Why a little bit goes a long way: Logical foundations of
  scientifically applicable mathematics.
\newblock In D.~Hull, M.~Forbes, and K.~Okruhlik, editors, \emph{{PSA}:
  Proceedings of the Biennial Meeting of the Philosophy of Science
  Association}, volume 1992, Volume Two: Symposia and Invited Papers, pages
  442--455. The University of Chicago Press on behalf of the Philosophy of
  Science Association, 1992.
\newblock \doi{10.1086/psaprocbienmeetp.1992.2.192856}.

\bibitem[Feferman(1998)]{Feferman1998}
S.~Feferman.
\newblock \emph{In the Light of Logic}.
\newblock Logic and Computation in Philosophy. Oxford University Press, New
  York and Oxford, 1998.

\bibitem[Feferman(2005)]{Feferman2005}
S.~Feferman.
\newblock Predicativity.
\newblock In S.~Shapiro, editor, \emph{The Oxford Handbook of Philosophy of
  Mathematics and Logic}, pages 590--624. Oxford University Press, Oxford,
  2005.
\newblock \doi{10.1093/0195148770.003.0019}.

\bibitem[Fenstad(1970)]{Fenstad1970}
J.~E. Fenstad, editor.
\newblock \emph{Selected Works Selected Works in Logic by Th.~Skolem}.
\newblock Universitetsforlaget, Oslo, Bergen, Troms{\"o}, 1970.

\bibitem[Field(1980/2016)]{Field2016}
H.~Field.
\newblock \emph{Science Without Numbers}.
\newblock Princeton University Press, 2nd edition, 1980/2016.

\bibitem[Field(1989)]{Field1989}
H.~Field.
\newblock \emph{Realism, Mathematics, and Modality}.
\newblock Blackwell, Oxford, 1989.

\bibitem[Field(1991)]{Field1991}
H.~Field.
\newblock Metalogic and modality.
\newblock \emph{Philosophical Studies}, 62\penalty0 (1):\penalty0 1--22, Apr
  1991.
\newblock \doi{10.1007/BF00646253}.

\bibitem[Fishburn(1970)]{Fishburn1970}
P.~C. Fishburn.
\newblock {Arrow}'s impossibility theorem: Concise proof and infinite voters.
\newblock \emph{Journal of Economic Theory}, 2\penalty0 (1):\penalty0 103--106,
  March 1970.
\newblock \doi{10.1016/0022-0531(70)90015-3}.

\bibitem[Friedman(1975)]{Friedman1975}
H.~Friedman.
\newblock Some systems of second order arithmetic and their use.
\newblock In R.~D. James, editor, \emph{Proceedings of the 17th International
  Congress of Mathematicians, Vancouver 1974}, volume~1, pages 235--242.
  International Mathematical Union, 1975.
\newblock URL \url{https://www.mathunion.org/fileadmin/ICM/Proceedings/
  ICM1974.1/ICM1974.1.ocr.pdf}.

\bibitem[Friedman(1976)]{Friedman1976}
H.~Friedman.
\newblock Systems of second order arithmetic with restricted induction. {I,
  II}.
\newblock \emph{The Journal of Symbolic Logic}, 41\penalty0 (2):\penalty0
  557--559, 1976.
\newblock \doi{10.1017/S0022481200051665}.
\newblock (Abstracts).

\bibitem[Ganea(2010)]{Ganea2010}
M.~Ganea.
\newblock Two (or three) notions of finitism.
\newblock \emph{The Review of Symbolic Logic}, 3\penalty0 (1):\penalty0
  119--144, 2010.
\newblock \doi{10.1017/S1755020309990323}.

\bibitem[Giaquinto(1983)]{Giaquinto1983}
M.~Giaquinto.
\newblock {Hilbert}'s philosophy of mathematics.
\newblock \emph{The British Journal for the Philosophy of Science}, 34\penalty0
  (2):\penalty0 119--132, 1983.
\newblock \doi{10.1093/bjps/34.2.119}.

\bibitem[Grzegorczyk(1955)]{Grzegorczyk1955}
A.~Grzegorczyk.
\newblock Elementarily definable analysis.
\newblock \emph{Fundamenta Mathematicae}, 41\penalty0 (2):\penalty0 311--338,
  1955.
\newblock \doi{10.4064/fm-41-2-311-338}.

\bibitem[Hallett(1979)]{Hallett1979}
M.~Hallett.
\newblock Towards a theory of mathematical research programmes {(I)}.
\newblock \emph{The British Journal for the Philosophy of Science}, 30\penalty0
  (1):\penalty0 1--25, 1979.
\newblock \doi{10.1093/bjps/30.1.1}.

\bibitem[Hallett(1984)]{Hallett1984}
M.~Hallett.
\newblock \emph{Cantorian Set Theory and Limitation of Size}.
\newblock Number~10 in Oxford Logic Guides. Oxford University Press, Oxford,
  1984.

\bibitem[Heath(1908)]{Heath1908}
T.~L. Heath.
\newblock \emph{The Thirteen Books of {Euclid}'s Elements}, volume~1.
\newblock Cambridge University Press, Cambridge, 1908.

\bibitem[Hellman(1989)]{Hellman1989}
G.~Hellman.
\newblock \emph{Mathematics Without Numbers: Towards a Modal-Structural
  Interpretation}.
\newblock Clarendon Press, Oxford, 1989.

\bibitem[Hellman(1999)]{Hellman1999}
G.~Hellman.
\newblock Some ins and outs of indispensability: A modal-structural
  perspective.
\newblock In A.~Cantini, E.~Casari, and P.~Minari, editors, \emph{Logic and
  Foundations of Mathematics: Selected Contributed Papers of the Tenth
  International Congress of Logic, Methodology and Philosophy of Science,
  Florence, August 1995}, Synthese Library, pages 25--39. Springer, Dordrecht,
  1999.
\newblock \doi{10.1007/978-94-017-2109-7}.

\bibitem[Hilbert and Bernays(1939)]{HilbertBernays1939}
D.~Hilbert and P.~Bernays.
\newblock \emph{Grundlagen der Mathematik}, volume~2.
\newblock Springer, Berlin, 1939.

\bibitem[Hilbert and Bernays(1970)]{HilbertBernays1970}
D.~Hilbert and P.~Bernays.
\newblock \emph{Grundlagen der Mathematik}, volume~II.
\newblock Springer-Verlag, Berlin, Heidelberg, New York, 2nd edition, 1970.

\bibitem[Hirschfeldt(2014)]{Hirschfeldt2014}
D.~R. Hirschfeldt.
\newblock \emph{{Slicing the Truth: On the Computability Theoretic and Reverse
  Mathematical Analysis of Combinatorial Principles}}.
\newblock Number~28 in Lecture Notes Series, Institute for Mathematical
  Sciences, National University of Singapore. World Scientific Publishing
  Company, 2014.
\newblock \doi{10.1142/9208}.

\bibitem[Irvine and Deutsch(2020)]{IrvineDeutsch2020}
A.~D. Irvine and H.~Deutsch.
\newblock Russell's paradox.
\newblock In E.~N. Zalta, editor, \emph{The {Stanford} Encyclopedia of
  Philosophy}. Metaphysics Research Lab, Stanford University, winter 2020
  edition, 2020.

\bibitem[Jockusch and Soare(1972)]{JockuschSoare1972}
C.~G. Jockusch and R.~I. Soare.
\newblock {$\Pi^0_1$} classes and degrees of theories.
\newblock \emph{Transactions of the American Mathematical Society},
  361:\penalty0 5805--5837, 1972.
\newblock \doi{10.1090/S0002-9947-1972-0316227-0}.

\bibitem[Kino et~al.(1970)Kino, Myhill, and Vesley]{KinoMyhillVesley1970}
A.~Kino, J.~Myhill, and R.~E. Vesley, editors.
\newblock \emph{Intuitionism and Proof Theory. Proceedings of the summer
  conference at Buffalo, N.Y., 1968}.
\newblock Number~60 in Studies in Logic and the Foundations of Mathematics.
  North-Holland, Amsterdam, 1970.

\bibitem[Kleene(1950)]{Kleene1950}
S.~C. Kleene.
\newblock A symmetric form of {G\"{o}del}'s theorem.
\newblock \emph{Koninklijke Nederlandse Akademie van Wetenschappen, Proceedings
  of the section of sciences}, 53:\penalty0 800--802, 1950.

\bibitem[Kleene(1952)]{Kleene1952a}
S.~C. Kleene.
\newblock Recursive functions and intuitionistic mathematics.
\newblock In L.~M. Graves, E.~Hille, P.~A. Smith, and O.~Zariski, editors,
  \emph{Proceedings of the International Congress of Mathematicians 1950},
  volume~II. American Mathematical Society, Providence, Rhode Island, 1952.

\bibitem[Kond\^{o}(1958)]{Kondo1958}
M.~Kond\^{o}.
\newblock {Sur les Ensembles Nommables et le Fondement de l'Analyse
  Math\'{e}matique, I}.
\newblock \emph{Japanese Journal of Mathematics}, 28:\penalty0 1--116, 1958.
\newblock \doi{10.4099/jjm1924.28.0_1}.

\bibitem[Kreisel(1958)]{Kreisel1958}
G.~Kreisel.
\newblock Mathematical significance of consistency proofs.
\newblock \emph{The Journal of Symbolic Logic}, 23\penalty0 (2):\penalty0
  155--182, June 1958.
\newblock \doi{10.2307/2964396}.

\bibitem[Kreisel(1959)]{Kreisel1959}
G.~Kreisel.
\newblock Analysis of the {Cantor--Bendixson} theorem by means of the analytic
  hierarchy.
\newblock \emph{Bulletin de l'Acad\'emie Polonaise des Sciences. S\'erie des
  Sciences Math\'ematiques, Astronomiques et Physiques}, 7:\penalty0 621--626,
  1959.

\bibitem[Kreisel(1960)]{Kreisel1960}
G.~Kreisel.
\newblock La predicativit{\'e}.
\newblock \emph{Bulletin de la Soci{\'e}t{\'e} Math{\'e}matiques de France},
  88:\penalty0 371--391, 1960.

\bibitem[Kreisel(1962)]{Kreisel1962}
G.~Kreisel.
\newblock The axiom of choice and the class of hyperarithmetic functions.
\newblock \emph{Indagationes Mathematicae (Proceedings)}, 65:\penalty0
  307--319, 1962.
\newblock \doi{10.1016/S1385-7258(62)50029-2}.

\bibitem[Kreisel(1967)]{Kreisel1967a}
G.~Kreisel.
\newblock Mathematical logic: What has it done for the philosophy of
  mathematics?
\newblock In R.~Schoenman, editor, \emph{Bertrand Russell, Philosopher of the
  Century}, pages 201--272. Allen and Unwin, 1967.

\bibitem[Kreisel(1968)]{Kreisel1968}
G.~Kreisel.
\newblock {A Survey of Proof Theory}.
\newblock \emph{The Journal of Symbolic Logic}, 33\penalty0 (3):\penalty0
  321--388, 1968.
\newblock \doi{10.2307/2270324}.

\bibitem[Kreisel(1970)]{Kreisel1970}
G.~Kreisel.
\newblock Church's thesis: a kind of reducibility axiom for constructive
  mathematics.
\newblock In A.~Kino, J.~Myhill, and R.~E. Vesley, editors, \emph{Intuitionism
  and Proof Theory. Proceedings of the summer conference at Buffalo, N.Y.,
  1968}, number~60 in Studies in Logic and the Foundations of Mathematics,
  pages 121--150. North-Holland, Amsterdam, 1970.
\newblock \doi{10.1016/S0049-237X(08)70746-8}.

\bibitem[Kreisel and Lacombe(1957)]{KreiselLacombe1957}
G.~Kreisel and D.~Lacombe.
\newblock Ensembles r{\'e}cursivement mesurables et ensembles r{\'e}cursivement
  ouverts ou ferm{\'e}s.
\newblock \emph{Comptes Rendus Hebdomadaires des S{\'e}ances de l'Acad{\'e}mie
  des Sciences}, 245\penalty0 (1):\penalty0 1106--1109, 1957.

\bibitem[{Le Roux} and Ziegler(2008)]{LeRouxZiegler2008}
S.~{Le Roux} and M.~Ziegler.
\newblock Singular coverings and non-uniform notions of closed set
  computability.
\newblock \emph{Mathematical Logic Quarterly}, 54\penalty0 (5):\penalty0
  545--560, September 2008.
\newblock \doi{10.1002/malq.200610058}.

\bibitem[Lewis(1920)]{Lewis1920}
F.~P. Lewis.
\newblock History of the parallel postulate.
\newblock \emph{The American Mathematical Monthly}, 27\penalty0 (1):\penalty0
  16--23, January 1920.
\newblock \doi{10.2307/2973238}.

\bibitem[Lorenzen(1955)]{Lorenzen1955}
P.~Lorenzen.
\newblock \emph{Einf{\"u}hrung in die Operative Logik und Mathematik}.
\newblock Springer, Berlin, Heidelberg, 1955.

\bibitem[Mancosu(1998)]{Mancosu1998}
P.~Mancosu, editor.
\newblock \emph{From Brouwer to Hilbert. The Debate on the Foundations of
  Mathematics in the 1920s}.
\newblock Oxford University Press, Oxford, 1998.

\bibitem[Montalb{\'a}n(2006)]{Montalban2006}
A.~Montalb{\'a}n.
\newblock Indecomposable linear orderings and hyperarithmetical analysis.
\newblock \emph{Journal of Mathematical Logic}, 6:\penalty0 89--120, 2006.
\newblock \doi{10.1142/S0219061306000517}.

\bibitem[Mostowski(1959)]{Mostowski1959}
A.~Mostowski.
\newblock On various degrees of constructivism.
\newblock In A.~Heyting, editor, \emph{Constructivity in Mathematics.
  Proceedings of the colloquium held at Amsterdam, 1957}, number 93, part B in
  Studies in Logic and the Foundations of Mathematics, pages 359--375.
  North-Holland, Amsterdam, 1959.
\newblock \doi{10.1016/S0049-237X(09)70270-8}.
\newblock Reprinted in {\citet[pp.~359--375]{Mostowski1979a}}.

\bibitem[Mostowski(1979)]{Mostowski1979a}
A.~Mostowski.
\newblock \emph{Foundational Studies: Selected Works. Volume 2}.
\newblock Number 93, Part B in Studies in Logic and the Foundations of
  Mathematics. North-Holland, Amsterdam, 1979.

\bibitem[Parsons(1970)]{Parsons1970}
C.~Parsons.
\newblock On a number-theoretic choice schema and its relation to induction.
\newblock In A.~Kino, J.~Myhill, and R.~E. Vesley, editors, \emph{Intuitionism
  and Proof Theory. Proceedings of the summer conference at Buffalo, N.Y.,
  1968}, number~60 in Studies in Logic and the Foundations of Mathematics,
  pages 459--473. North-Holland, Amsterdam, 1970.

\bibitem[Poincar\'{e}(1902)]{Poincare1902}
H.~Poincar\'{e}.
\newblock \emph{La Science et L'Hypoth\`{e}se}.
\newblock Flammarion, 1902.
\newblock English translation in {\citet{Poincare1913}}.

\bibitem[Poincar\'{e}(1905)]{Poincare1905}
H.~Poincar\'{e}.
\newblock \emph{Le Valeur de la Science}.
\newblock Flammarion, 1905.
\newblock English translation in {\citet{Poincare1913}}.

\bibitem[Poincar\'{e}(1908)]{Poincare1908}
H.~Poincar\'{e}.
\newblock \emph{Science et Methode}.
\newblock Flammarion, 1908.
\newblock English translation in {\citet{Poincare1913}}.

\bibitem[Poincar{\'e}(1910)]{Poincare1910}
H.~Poincar{\'e}.
\newblock {{\"U}ber transfinite Zahlen}.
\newblock In \emph{Sechs Vortr{\"a}ge {\"u}ber ausgew{\"a}hlte Gegenst{\"a}nde
  aus der reinen Mathematik und mathematischen Physik}, pages 45--48. Teubner,
  Leipzig/Berlin, 1910.
\newblock English translation in \citet[pp.~1071--1074]{Ewald1996b}.

\bibitem[Poincar\'{e}(1913)]{Poincare1913}
H.~Poincar\'{e}.
\newblock \emph{The Foundations of Science: The Foundations of Science: Science
  and Hypothesis, The Value of Science, Science and Method}.
\newblock Number~1 in Science and Education. The Science Press, 1913.
\newblock English translation by G.~B.~Halstead of {\citet{Poincare1902,
  Poincare1905, Poincare1908}}.

\bibitem[Pour-El and Richards(1989)]{Pour-ElRichards1989}
M.~Pour-El and J.~Richards.
\newblock \emph{Computability in Analysis and Physics}.
\newblock Perspectives in Mathematical Logic. Springer-Verlag, 1989.
\newblock xi + 206 pages.

\bibitem[Putnam(1971)]{Putnam1971}
H.~Putnam.
\newblock \emph{Philosophy of Logic}.
\newblock Harper, New York, 1971.

\bibitem[Quine(1970/1986)]{Quine1986}
W.~V. Quine.
\newblock \emph{Philosophy of Logic}.
\newblock Harvard University Press, second edition, 1970/1986.

\bibitem[Quine(1998)]{Quine1998}
W.~V. Quine.
\newblock Reply to {Parsons}.
\newblock In L.~E. Hahn and P.~A. Schilpp, editors, \emph{The Philosophy of W.
  V. Quine}, number~18 in Library of Living Philosophers, pages 398--403. Open
  Court, 2nd expanded edition, 1998.

\bibitem[Riker(1982)]{Riker1982}
W.~H. Riker.
\newblock \emph{Liberalism Against Populism: A Confrontation Between the Theory
  of Democracy and the Theory of Social Choice}.
\newblock Waveland Press, Long Grove, IL, 1982.

\bibitem[Russell(1907)]{Russell1907a}
B.~Russell.
\newblock On some difficulties in the theory of transfinite numbers and order
  types.
\newblock \emph{Proceedings of the London Mathematical Society, Series 2},
  4\penalty0 (1):\penalty0 29--53, 1907.
\newblock \doi{10.1112/plms/s2-4.1.29Digital Object Identifier (DOI)}.

\bibitem[Russell(1908)]{Russell1908}
B.~Russell.
\newblock Mathematical logic as based on the theory of types.
\newblock \emph{American Journal of Mathematics}, 30\penalty0 (3):\penalty0
  222--262, July 1908.
\newblock \doi{10.2307/2369948}.

\bibitem[Sacks(1990)]{Sacks1990}
G.~Sacks.
\newblock \emph{Higher Recursion Theory}.
\newblock Perspectives in Mathematical Logic. Springer-Verlag, 1990.

\bibitem[Shapiro(1991)]{Shapiro1991}
S.~Shapiro.
\newblock \emph{Foundations Without Foundationalism: A Case for Second-Order
  Logic}.
\newblock Number~17 in Oxford Logic Guides. Oxford University Press, Oxford,
  1991.
\newblock \doi{10.1093/0198250290.001.0001}.

\bibitem[Shore(2010)]{Shore2010}
R.~A. Shore.
\newblock Reverse mathematics: The playground of logic.
\newblock \emph{The Bulletin of Symbolic Logic}, 16\penalty0 (3):\penalty0
  378--402, 2010.
\newblock \doi{10.2178/bsl/1286284559}.

\bibitem[Sieg(1985)]{Sieg1985a}
W.~Sieg.
\newblock Reductions of theories for analysis.
\newblock In G.~Dorn and P.~Weingartner, editors, \emph{Foundations of Logic
  and Linguistics: Problems and Their Solutions}, pages 199--231.
  Springer-Verlag, Boston, MA, 1985.
\newblock \doi{10.1007/978-1-4899-0548-2_9}.

\bibitem[Simpson(1985)]{Simpson1985}
S.~G. Simpson.
\newblock {Friedman}'s research on subsystems of second order arithmetic.
\newblock In L.~A. Harrington, M.~D. Morley, A.~Scedrov, and S.~G. Simpson,
  editors, \emph{Harvey Friedman's Research on the Foundations of Mathematics},
  number 117 in Studies in Logic and the Foundations of Mathematics, pages
  137--159. North-Holland, Amsterdam, 1985.
\newblock \doi{10.1016/S0049-237X(09)70158-2}.

\bibitem[Simpson(1988)]{Simpson1988}
S.~G. Simpson.
\newblock {Partial realizations of Hilbert's program}.
\newblock \emph{The Journal of Symbolic Logic}, 53:\penalty0 349--363, 1988.
\newblock \doi{10.1017/S0022481200028309}.

\bibitem[Simpson(2009)]{Simpson2009}
S.~G. Simpson.
\newblock \emph{Subsystems of Second Order Arithmetic}.
\newblock Perspectives in Logic. Association for Symbolic Logic and Cambridge
  University Press, Cambridge, 2nd edition, 2009.
\newblock \doi{10.1017/cbo9780511581007}.

\bibitem[Skolem(1920)]{Skolem1920}
T.~Skolem.
\newblock Logisch-kombinatorische untersuchungen {\"u}iber die
  erf{\"u}llbarkeit oder beweisbarkeit mathematischer satze nebst einem
  theoreme {\"u}ber dichte mengen.
\newblock \emph{Skrifter utgitt av Videnskapsselskapet i Kristiana, I, Math.
  Naturv. KL}, 4, 1920.

\bibitem[Skolem(1923)]{Skolem1923a}
T.~Skolem.
\newblock {Einige Bemerkungen zur axiomatischen Begr{\"u}ndung der
  Mengenlehre}.
\newblock In \emph{Matematikerkongressen i Helsingfors 4--7 Juli 1922. Den
  femte skandinaviska matematikerkongressen, Redog{\"o}relse}, pages 217--232.
  Akademiska Bokhandeln, Helsingfors, 1923.
\newblock Reprinted in {\citet[pp.~137--152]{Fenstad1970}}. Translated into
  English as `Some remarks on axiomatized set theory'' by S.~Bauer-Mengelberg
  in {\citep[pp.~290--301]{vanHeijenoort1967}}.

\bibitem[Soare(2016)]{Soare2016}
R.~I. Soare.
\newblock \emph{Turing Computability}.
\newblock Theory and Applications of Computability. Springer, Berlin,
  Heidelberg, 2016.
\newblock \doi{10.1007/978-3-642-31933-4}.

\bibitem[Specker(1949)]{Specker1949}
E.~Specker.
\newblock {Nicht konstruktiv beweisbare S\"{a}tze der Analysis}.
\newblock \emph{The Journal of Symbolic Logic}, 14\penalty0 (3):\penalty0
  145--158, 1949.
\newblock \doi{10.2307/2267043}.

\bibitem[Tait(1981)]{Tait1981}
W.~W. Tait.
\newblock Finitism.
\newblock \emph{The Journal of Philosophy}, 78\penalty0 (9):\penalty0 524--546,
  1981.
\newblock \doi{10.2307/2026089}.

\bibitem[{van Atten}(2004)]{vanAtten2004}
M.~{van Atten}.
\newblock \emph{On Brouwer}.
\newblock Wadsworth Philosophers Series. Wadsworth/Thomson Learning, Belmont,
  CA, 2004.

\bibitem[{van Heijenoort}(1967)]{vanHeijenoort1967}
J.~{van Heijenoort}, editor.
\newblock \emph{From Frege to G{\"o}del: A Source Book in Mathematical Logic,
  1879--1931}.
\newblock Harvard University Press, Cambridge, MA, 1967.

\bibitem[Wang(1954)]{Wang1954}
H.~Wang.
\newblock The formalization of mathematics.
\newblock \emph{Journal of Symbolic Logic}, 19\penalty0 (4):\penalty0 241--266,
  December 1954.
\newblock \doi{10.2307/2267732}.

\bibitem[Weaver(2009)]{Weaver2009}
N.~Weaver.
\newblock {Predicativity beyond {$\Gamma_0$}}.
\newblock \href{http://arxiv.org/abs/math/0509244v3}{arXiv:math/0509244v3},
  2009.

\bibitem[Weber and Colyvan(2010)]{WeberColyvan2010}
Z.~Weber and M.~Colyvan.
\newblock A topological sorites.
\newblock \emph{The Journal of Philosophy}, 107\penalty0 (6):\penalty0
  311--325, 2010.
\newblock \doi{10.5840/jphil2010107624}.

\bibitem[Weyl(1918)]{Weyl1918}
H.~Weyl.
\newblock \emph{Das Kontinuum: Kritische Untersuchungen {\"u}ber die Grundlagen
  der Analysis}.
\newblock Veit, 1918.
\newblock Reprinted in: H.\ Weyl, E.\ Landau, and B.\ Riemann, \emph{Das
  Kontinuum und andere Monographien}, Chelsea, 1960, 1973.

\bibitem[Weyl(1925)]{Weyl1925}
H.~Weyl.
\newblock {Die heutige Erkenntnislage in der Mathematik}.
\newblock \emph{Symposion}, 1:\penalty0 1--32, 1925.
\newblock Reprinted in {\citep[p.~511--42]{Weyl1968}}. English translation in
  {\citep[pp.~123--42]{Mancosu1998}}.

\bibitem[Weyl(1968)]{Weyl1968}
H.~Weyl.
\newblock \emph{Gesammelte Abhandlungen}.
\newblock Springer, Berlin and Heidelberg, 1968.
\newblock 4 vols.

\bibitem[Zermelo(1904)]{Zermelo1904}
E.~Zermelo.
\newblock {Neuer Beweis, da\ss{} jede Menge Wohlordnung werden kann (Aus einem
  an Herrn Hilbert gerichteten Briefe)}.
\newblock \emph{Mathematische Annalen}, 59:\penalty0 514--516, December 1904.
\newblock \doi{10.1007/BF01445300}.
\newblock {English} translation by {Stefan Bauer-Mengelberg} in
  {\citep[pp.~139--141]{vanHeijenoort1967}}.

\bibitem[Zermelo(1908)]{Zermelo1908}
E.~Zermelo.
\newblock {Neuer Beweis f\"{u}r die M\"{o}glichkeit einer Wohlordnung}.
\newblock \emph{Mathematische Annalen}, 65\penalty0 (1):\penalty0 107--128,
  1908.
\newblock \doi{10.1007/BF01450054}.
\newblock English translation by Stefan Bauer-Mengelberg in
  {\citep[pp.~183--198]{vanHeijenoort1967}}.

\end{thebibliography}

\end{document}